\documentclass[12pt]{extarticle}
\usepackage{amsmath, amsthm,graphicx,amsfonts,amssymb,mathrsfs,showidx}

\usepackage{relsize}
\usepackage{color}
\definecolor{rojo}{rgb}{1,0,0}
\usepackage[all]{xy}
\makeindex
\usepackage{hyperref}
\usepackage{color}
\definecolor{abelian}{cmyk}{0.50,0,1,.4}
\definecolor{noabelian}{cmyk}{0.94,0.54,0,0}
\definecolor{rojo}{cmyk}{0,1,1,0}
\definecolor{verde}{cmyk}{0.91,0,0.88,0.12}
\xyoption{arc}

\usepackage{tikz}
\usetikzlibrary{tqft}

\newtheorem{thm}{Theorem}

\newtheorem{prop}[thm]{Proposition}

\theoremstyle{definition}
\newtheorem{defn}[thm]{\textbf{Definition}}

\theoremstyle{definition}

\theoremstyle{remark}
\newtheorem{rem}[thm]{Remark:}

\theoremstyle{Notation}

\newcommand{\op}{\operatorname}

\newcommand{\ben}{\begin{equation}}
\newcommand{\een}{\end{equation}}
\newcommand{\bena}{\begin{equation*}}
\newcommand{\eena}{\end{equation*}}

\newcommand{\ma}{\mathcal}
\newcommand{\To}{\longrightarrow}

\newcommand{\cob}{\mathscr{S}}
\newcommand{\cobG}{\mathscr{S}^G}
\newcommand{\Ini}{\tilde{\cob}_0^G}
\newcommand{\Gru}{\mathscr{G}^G}

\newcommand{\ZZ}{\mathbb{Z}}
\newcommand{\Tot}{\longmapsto}
\def\NN{\mathbb{N}}

\def\ZZ{\mathbb{Z}}
\def\CC{\mathbb{C}}

\title{The classifying space of the 1+1 dimensional \\ $G$-cobordism category}
\author{Carlos Segovia\footnote{
Instituto de Matem\'aticas UNAM-Oaxaca,
\newline
csegovia@matem.unam.mx}}
\begin{document}
\maketitle

\begin{abstract}
For a finite group $G$, we define the $G$-cobordism category in dimension two.
We show there is a one-to-one correspondence between the connected components of its classifying space and the abelianization of $G$. 
Also, we find an isomorphism of its fundamental group onto the direct sum  $\ZZ\oplus \Omega_2^{SO}(BG)$, where $\Omega_2^{SO}(BG)$ is the 
free oriented $G$-bordism group in dimension two, 
and we study the classifying space of some important subcategories.
We obtain the classifying space 
has the homotopy type of the product $G/[G,G]\times S^1\times X^G$, where 
$\pi_1(X^G)=\Omega_2^{SO}(BG)$.
Finally, we present some results about the classification of $G$-topological quantum field theories in dimension two.

\end{abstract}

\section*{Introduction}
\label{intro}\label{sec1}

Let $\cobG$ denote the 1+1 dimensional $G$-cobordism category. The objects of $\cobG$ are
disjoint unions of principal $G$-bundles over the circle.
They are represented by finite sequences $(x_1,x_2,\cdots,x_n)$, $x_i\in G\sqcup \{0\}$, where $x_i\in G$ denotes the monodromy of a principal $G$-bundle over the circle and $x_i=0$ denotes the empty $G$-bundle. When $x_i$ is equal to the neutral $1\in G$, we say that we have trivial monodromy.
The morphisms of $\cobG$ are cobordisms 
of principal $G$-bundles over oriented surfaces up to diffeomorphism.

Composition of two cobordisms is given by the gluing of their collars for a common boundary. The gluing of two manifolds has a smooth structure well defined up to diffeomorphism, which extends to the gluing of two principal $G$-bundles producing the appropriate cobordism. 

Associated to $\cobG$ we have a simplicial set $N\cobG$, called its nerve \cite{segal}. The geo\-me\-tric realization of $N \cobG$ is denoted by $B\cobG$. 
We show in Section \ref{seccob} that the connected components of $B\cobG$ have the form described in the following theorem. 
\begin{thm}\label{thm1}There is a bijection $\pi_0(B\cobG)\cong G/[G,G]$ which is compatible with the tensor product of $\cobG$.
\end{thm}

Let $\Omega_2^{SO}(BG)$ be the free oriented G-bordism group in dimension two from Conner-Floyd \cite{CF}.
In Section \ref{sec3} we show that the fundamental group of $B\cobG$ is isomorphic with the direct sum.

\begin{thm}\label{confund1} There is an isomorphism $\pi_1(B\cobG)\cong \ZZ\oplus \Omega_2^{SO}(BG)$.
\end{thm}
The $G$-equivariant version of the embedded cobordism category of \cite{a6}, denoted by $\mathcal{C}_2^G$, has connected components $\pi_0(B\mathcal{C}_2^G)\cong \Omega_1^{SO}(BG)=G/[G,G]$. B\"okstedt-Svane \cite{a666} show that 
the fundamental group $\pi_1(B\mathcal{C}_2^G)$ is obtained in terms of generators up to the ``chimera relations". In Section \ref{sec3} we show these relations are equivalent to the generating relations for the category of fractions of $\cobG$. 
Thus $\mathcal{C}_2^G$ and $\cobG$ have isomorphic fundamental groups. 

For $G$ the trivial group, Tillmann \cite{a3} proves that $B\cob\simeq S^1\times X$ where $X$ is a simply connected infinite loop space. As a consequence, the homotopy groups satisfy $\pi_n(B\cob)\cong \pi_n(X)$ for $n\geq 2$. Juer-Tillman conjecture that indeed $\pi_n(X)=0$, for $n\geq 2$, as a consequence of the discrete localisation conjecture \cite[p.217]{loca}. 
More precisely, this is equivalent to show that the canonical inclusion $\cob_b\hookrightarrow \cob$ is a homotopy equivalence, where $\cob_b$ is the positive boundary subcategory ($B\cobG_b\simeq S^1$). In fact, Tillmann \cite{a3} shows that the composition $\cob_b\hookrightarrow\cob\longrightarrow\cob[\cob^{-1}]$ is a homotopy equivalence.
In the embedded cobordism category \cite[Thm 6.1]{a6} this inclusion is a weak equivalence.
We do not know the behaviour for the other homotopy groups, for instance, we ignore if the canonical map between these two categories induces an epimorphism for the two dimensional homotopy groups.
We have the same question for the $G$-equivariant case.

The $G$-cobordism category $\cobG$ is a symmetric monoidal category, hence the classifying space $B\cobG$ is an $H$-space. Theorem \ref{thm1} implies that the connected components has a group structure, hence every pair of connected components are homotopy equivalent by multiplication by an element. Moreover, Theorem \ref{confund1} implies that $B\cobG$ has abelian fundamental group.
Thus $B\cobG$ is homotopic to an infinite loop space, which splits into the product of $G/[G,G]$ and some connected component of $B\cobG$. 
However, we need to restrict the category $\cobG$ to certain type of morphisms in order to extend the Tillmann's procedure \cite{a3} (this is specified in Section \ref{sec2}). We show in Section \ref{sec4} that the classifying space of the $G$-cobordism category has the homotopy type $G/[G,G]\times S^1\times X^G$, where $\pi_1(X^G)=\Omega_2^{SO}(BG)$.
In this case the discrete localisation conjecture has the form:
\begin{itemize}
\item The homotopy type of the $G$-cobordism category is $G/[G,G]\times S^1\times X^G$, where $X^G$ is a $K(G,1)$.
\end{itemize}

There are a number of other papers which study bordism categories of surfaces with principal $G$-bundles. The works that study the same category as in our paper are the dissertation by Sweet \cite{sweet}, the papers by Kaufmann \cite{kaufmann,Kauffmann2}, the work of Moore-Segal \cite{a6} and the paper by Gonz\'alez-Segovia \cite{GS}. 
There is a different perspective which results in an equivalent bordism category by Turaev \cite{turaev} and Tagami \cite{tagami}. This point of view considers objects and morphisms provided with a map to the Eilenberg-MacLane space $K(G,1)$. Lastly, the bordism category studied by Davidovich \cite{davi} is different from our bordism category since it comes in the context of $(\infty,n)$-category for the cobordism hypothesis.

This article is organized as follows:
in Section \ref{seccob} we define the $G$-cobordism category $\cobG$. 
In Section \ref{sec3} we describe the groupoid of fractions of $\cobG$ and we find the fundamental group of $\cobG$. 
In Section \ref{sec2} we study the classifying spaces of some important 
subcategories of $\cobG$. 
%
In Section \ref{sec4} we study the homotopy type of the classifying space of the $G$-cobordism category. 
Finally, in Section \ref{sec5} we give some results about the classification of $G$-topological quantum field theories.

\section{The $G$-cobordism category}
\label{seccob}

Every principal $G$-bundle over the circle is isomorphic to the following construction:
attach the ends of $[0,1]\times G$ via multiplication by some element $x\in G$, i.e., $(0,y)$ is identified with $(1,yx)$ for every $y\in G$. 
This construction projects to the circle by restricting to the first coordinate and the action of $G$ is defined by left multiplication  on the second coordinate.
 Throughout the paper the element $x\in G$ is called the {\it monodromy} of the corresponding principal $G$-bundle over the circle. When $x$ is the neutral element, $1\in G$, we say that we have trivial monodromy.
Additionally, we consider the empty $G$-bundle represented by the zero $0$.

 A finite sequence $(x_1,x_2,\cdots,x_n)$, with $x_i\in G\sqcup \{0\}$, $1\leq i\leq n$, represents the disjoint union of principal $G$-bundles
 over the circle for $x_i\in G$ and the empty $G$-bundle for $x_i=0$. In the case $n=1$, we can take out the parenthesis.

\begin{figure}
    \centering
    \includegraphics{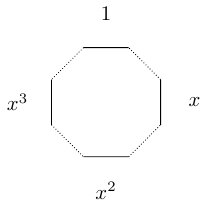}
    \hspace{1cm}
\includegraphics{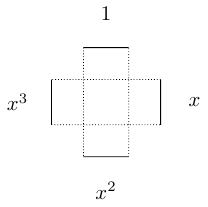}    
\caption{Principal $\ZZ_4$-bundles with monodromies $x$ and $x^2$, respectively.}
    \label{f2}
\end{figure}

\begin{figure}
    \centering
    \includegraphics{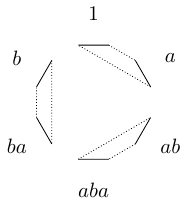} 
    \hspace{1cm}
     \includegraphics{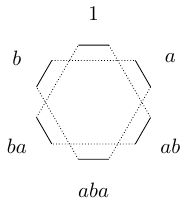}
    \caption{Principal $S_3$-bundles with monodromies $a$ and $ab$, respectively.}
    \label{f23}
\end{figure}

For example, we consider the cyclic group $\ZZ_4=\{1,x,x^2,x^3\}$ and the symmetric group $S_3=\langle a,b:a^2=b^2=(ab)^3=1\rangle$. For $x,x^2\in\ZZ_4$ and $a,ab\in S_3$, we obtain the four principal $G$-bundles of Figures \ref{f2} and \ref{f23} (the dotted lines denote the identifications).

Let $\Sigma$ and $\Sigma'$ be $d$-dimensional closed smooth ma\-nifolds. A cobordism between $\Sigma$ and $\Sigma'$ is a $(d+1)$-dimensional smooth manifold $M$, with boundary diffeomorphic to $\Sigma\sqcup-\Sigma'$, where $-\Sigma'$ is $\Sigma'$ with the reverse orientation. Two cobordisms $M$ and $M'$ are equivalent if there exists a diffeomorphism $\phi:M\To M'$ such that we have the commutative diagram
\ben
\xymatrix{&M\ar[dd]^\phi&\\\Sigma\ar[ru]\ar[rd]&&\Sigma'\ar[lu]\ar[ld]\\&M'&\,.}
\een 

Let $\xi:P\rightarrow \Sigma$ and $\xi':P'\rightarrow \Sigma'$ be principal $G$-bundles. 
A {\it $G$-cobordism} between $\xi$ and $\xi'$ is a principal $G$-bundle $\epsilon:Q\rightarrow M$, with diffeomorphisms for the boundaries $\partial M\cong S\sqcup- S'$ and $\partial Q\cong P\sqcup -P'$, which match with the projections and the restriction of the action. 
Two $G$-cobordisms $\epsilon:Q\rightarrow M$ and ${\epsilon}':Q'\rightarrow M'$ define the same class if $M$ and ${M}'$ are equivalent as cobordisms by a diffeomorphism $\phi:M\rightarrow {M}'$, $Q$ and $Q'$ are equivalent as cobordisms by a $G$-equivariant diffeomorphism $\psi:Q\To Q'$, and in addition, we have the commutative diagram  
\ben
\xymatrix{Q\ar[r]^\psi\ar[d]_\epsilon&Q'\ar[d]^{\epsilon'}\\M\ar[r]_\phi & M'\,.}
\een

In dimension two we give some important examples of $G$-cobordisms:

\begin{enumerate}
\item (Cylinder) for $x,y\in G$, this is the $G$-cobordism from $x$ to $y$ with base space the cylinder, where there is an element $z\in G$ with $y=zxz^{-1}$. The diffeomorphism identification is given by the Dehn twist. Notice that the conjugation by the neutral $1\in G$ produces the identity morphism $\op{id}_x$. 

\item (Pair of pants) for $x,y\in G$, consider the principal $G$-bundle over the wedge of two circles with monodromies $(x,y)$. The pullback by the retraction of the pair of pants to $S^1\vee S^1$, defines a $G$-cobordism over the pair of pants with incoming boundary $(x,y)$ and outgoing boundary the product $xy$. We denote  this element by $P_{(x,y)}: (x,y)\rightarrow xy$. Similarly, we define the reverse $G$-cobordism over the pair of pants $\overline{P}_{(x,y)}:xy\rightarrow (x,y)$, which is $P_{(x,y)}$ with the opposite orientation.

\item (Pair of pants with multiple legs) for $x_i\in G$, with $1\leq i\leq n$, consider the principal $G$-bundle over the wedge $S^1\vee \stackrel{n}{\cdots}\vee S^1$ with monodromies $\hat{x}=(x_1,\cdots,x_n)$. The pullback by the retraction of the pair of pants with $n$-legs to 
$S^1\vee \stackrel{n}{\cdots}\vee S^1$, defines a $G$-cobordism over the pair of pants with $n$-legs with incoming boundary $\hat{x}=(x_1,\cdots,x_n)$ and outgoing boundary the product $\prod_{i=1}^n x_i$.
We represent this element by $P_{\hat{x}}$ and similarly, we define the $G$-cobordism $\overline{P}_{\hat{x}}$.

\item (Disk) there is only one $G$-cobordism over the disk which is a trivial bundle. We represent this element by $D:1\rightarrow 0$ and similarly, we define the $G$-cobordism $\overline{D}:0\rightarrow 1$.

\item (Handlebody) for $x_i,y_i\in G$, $1\leq i\leq n$, this is the $G$-cobordism with base space the two dimensional handlebody of genus $n$ with one boundary circle. The base space is constructed from a polygon with $4n$ sides where a disc is removed inside. 
The $4n$ sides have associated monodromies $y_n$, $x_n$, $y_n^{-1}$, $x_n^{-1}$, $\cdots$, $y_1$, $x_1$, $y_1^{-1}$, $x_1^{-1}$ in clockwise order. 
Therefore, the monodromy for the boundary circle is given by the product $[y_n,x_n]\cdots [y_1,x_1]$. 
\end{enumerate}

\begin{figure}[h!]
\begin{center}
\includegraphics{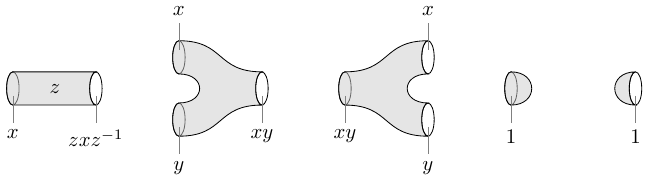}
 \end{center}
 \caption{$G$-cobordisms over the cylinder, the pair of pants and the disc.}
 \label{fig3}
 \end{figure}

\begin{figure}[h!]
\begin{center}
\includegraphics{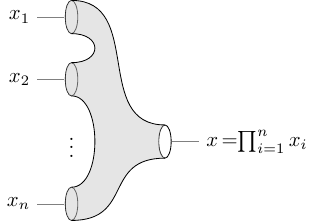}
 \end{center}
 \caption{$G$-cobordisms over the pair of pants with multiple legs $P_{(x_1,\cdots,x_n)}$.}
 \label{fig4}
 \end{figure}

\begin{figure}[h!]
\begin{center}
   \includegraphics{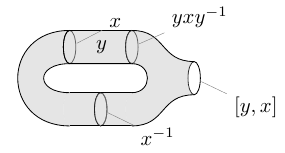}
 \end{center}
 \caption{$G$-cobordism over a handlebody with monodromy $[y,x]$ for the boundary circle.}
 \label{fig5}
 \end{figure}

In Figure \ref{fig3}, Figure \ref{fig4} and Figure \ref{fig5} we have pictures of the previous $G$-cobordisms. 
For these pictures, we draw from left to right the direction of our cobordisms. Also every circle is labelled with its monodromy and for the cylinders we include inside the element of the group with which we do the conjugation.

\begin{defn} \label{cobordismoe}The objects of the $G$-cobordism category $\cobG$ are finite sequences $(x_1,\cdots,x_n)$, $x_i\in G\sqcup\{0\}$,
where $x_i\in G$ denotes the monodromy of a principal $G$-bundle over the circle and $x_i=0$ denotes the empty $G$-bundle. The morphisms of $\cobG$ are $G$-cobordisms over oriented surfaces.
The morphisms in $\cobG$ have a well defined composition by the gluing of their collars for a common boundary.
\end{defn}


The category $\cobG$ has a monoidal structure given by the disjoint union $\sqcup: \cobG\times \cobG\rightarrow \cobG$, where the unit is the empty $G$-bundle. Furthermore, there is a symmetric structure induced by a natural transformation $c:\sqcup \rightarrow \sqcup \circ \tau$ where $\tau$ is the twist functor. This is given by a bunch of crossing identity $G$-cobordisms over  cylinders, where $c^2=\op{id}$.



Now, we show Theorem \ref{thm1} which says that there is one-to-one correspondence between the connected components of $B\cobG$ and the  abelianization of $G$, i.e,
\ben
\pi_0(B\cobG)\cong \frac{G}{[G,G]}\,.
\een
Because of the $G$-cobordism over the disc, we can restrict to objects given by finite sequences $(x_1,\cdots,x_n)$, with $x_i\neq 0$. Then the $G$-cobordism over the pair of pants with multiple legs connects $(x_1,\cdots,x_n)$ with its product $x=\prod_{i=1}^nx_i$. Due to the $G$-cobordisms over the handlebodies, two elements $x, y\in G$ are connected in $\cobG$ if and only if they differ by an element in $[G,G]$. Notice that the product in $\pi_0(B\cobG)$ induced by the monoidal structure of $\cobG$, agrees with the product of $G/[G,G]$. This is because of the existence of the $G$-cobordism over the pair of pants (with multiple legs).

\section{The fundamental group}
\label{sec3}


The aim of the present section is to show that the fundamental group of $B\cobG$ is the direct sum $\ZZ\oplus \Omega_2^{SO}(BG)$, where $\Omega_2^{SO}(BG)$ is the free oriented $G$-bordism group in dimension two. In this respect, we need an interpretation of the groupoid of fractions $\Gru:=\cobG[(\cobG)^{-1}]$ with the left calculus of fractions \cite{a9}. The fundamental group with base space an object $x$ in $\cobG$, is the automorphism group $\op{Hom}_{\Gru}(x,x)$, see Quillen \cite{a7}.

An important assumption is the restriction of $\cobG$ to objects which are in the connected component of the empty $G$-bundle.
Thus, for every object $\hat{x}$, we can take a morphism $\delta_{\hat{x}}$ from $\hat{x}$ to the empty $G$-bundle, where for simplicity, take it with connected base space.

Consider the subset of morphisms of those endomorphisms in $\cobG$ which are the disjoint union of identity $G$-cobordisms over cylinders and $G$-cobordisms over closed surfaces. Denote this subset of morphisms by $\Ini$.
Notice $\Ini$ admits a calculus of left fractions \cite{a9}. 
Therefore, the category of fractions $\cobG[ (\Ini)^{-1} ]$ can be described by roofs 
\ben\label{roofa1}
\xymatrix{(\op{id}_{\hat{x}}\sqcup\,\gamma,\Sigma)\hspace{0.5cm}=&&\hat{x}\ar[ld]_{\op{id}_{\hat{x}}\sqcup\,\gamma}\ar[rd]^\Sigma&\\ &\hat{x}&&\hat{y}\,,}
\een
where $\Sigma:\hat{x}\rightarrow\hat{y}$ is a morphism in $\cobG$ and $\gamma$ is a $G$-cobordism over closed surfaces. 

However, the category $\cobG[ (\Ini)^{-1} ]$ is not enough to obtain the groupoid of fractions $\Gru$. For this purpose, we consider the quotient by an equivalence relation which is motivated as follows: the functor $\cobG\rightarrow\cobG[ (\Ini)^{-1} ]$ sends a morphism  $\Sigma:\hat{x}\longrightarrow\hat{y}$ to the roof 
$(\op{id}_{\hat{x}},\Sigma)$. 
Take the gluing along $\hat{y}$ of $\Sigma$ and $\delta_{\hat{y}}$. This is represented by $\Sigma'=\Sigma\sqcup_{\hat{y}}\delta_{\hat{y}}$ which is a morphism from $\hat{x}$ to the empty $G$-bundle.
Denote by $\overline{\Sigma'}:0\rightarrow \hat{x}$ the $G$-cobordism obtained by changing the orientation of $\Sigma'$.
We can obtain a left inverse for $(\op{id}_{\hat{x}},\Sigma)$ after we impose some identification. The proposal for the left inverse is the roof $\left(\op{id}_{\hat{y}}\sqcup (\Sigma'\circ\overline{\Sigma'}), \overline{\Sigma'}\circ\delta_{\hat{y}}\right)$. In fact, the composition of 
$(\op{id}_{\hat{x}},\Sigma)$ followed by    
$\left(\op{id}_{\hat{y}}\sqcup (\Sigma'\circ\overline{\Sigma'}),\overline{\Sigma'}\circ\delta_{\hat{y}}\right)$ is 
\ben\left(\op{id}_{\hat{x}}\sqcup (\Sigma'\circ\overline{\Sigma'}),\overline{\Sigma'}\circ\Sigma'\right)\,.\een
Therefore, we need to impose the identification \ben\label{identification}\op{id}_{\hat{x}}\sqcup (\Sigma'\circ\overline{\Sigma'})\sim\overline{\Sigma'}\circ\Sigma'\,.\een
Take the equivalence relation generated by the identifications \eqref{identification} using the composition and the monoidal structure of $\cobG$ and denote the quotient category by $\cobG_\sim$.

\begin{figure}
    \centering
   \includegraphics{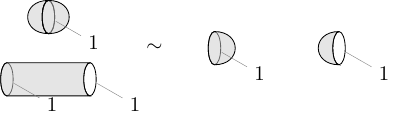}
       \caption{The identification $\op{id}_1\sqcup (D\circ\overline{D} )\sim \overline{D}\circ D$.}
    \label{higuro}
\end{figure}

\begin{figure}
    \centering
   \includegraphics{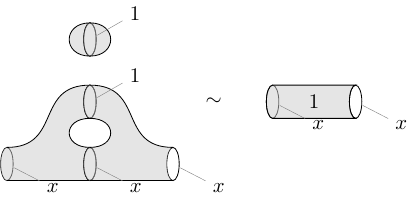}
    \caption{Property in the category of fractions.}
    \label{cancelito}
\end{figure}

In Figure \ref{higuro} we illustrate \eqref{identification} for the disc $D:1\rightarrow 0$ and the reverse disc $\overline{D}:0\rightarrow 1$. 
As an application we obtain that in the category of fractions we have some sort of cancellation between $G$-cobordisms over handlebodies and spheres, see Figure \ref{cancelito}.

In the following proposition we outline the missing details to calculate the category of fractions.
\begin{prop}\label{dfr} There is an isomorphism $\Gru\cong\cobG_\sim[(\Ini)^{-1}]$.
\end{prop}
\begin{proof}
Now we show that the left inverse of $(\op{id}_{\hat{x}},\Sigma)$ is also right inverse. For this purpose, we consider the composition of  
$(\op{id}_{\hat{y}}\sqcup (\Sigma'\circ\overline{\Sigma'}),\overline{\Sigma'}\circ\delta_{\hat{y}})$ followed by 
$(\op{id}_{\hat{x}},\Sigma)$. This is equal to  
\ben(\alpha,\beta):=\left(\op{id}_{\hat{y}}\sqcup(\Sigma'\circ\overline{\Sigma'}), (\Sigma\circ\overline{\Sigma})\sqcup(\delta_{\hat{y}}\circ\overline{\delta_{\hat{y}}})\right)\,.\een 
This is equivalent to the identity morphism because we can take a third roof  
\ben(\epsilon,\xi):=\left(\op{id}_{\hat{y}}\sqcup(\delta_{\hat{y}}\circ\overline{\delta_{\hat{y}}}),
(\Sigma\circ\overline{\Sigma})\sqcup(\delta_{\hat{y}}\circ\overline{\delta_{\hat{y}}})^{\sqcup 2}
\right)\,,\een and we have the commutative diagram 
\ben
\xymatrix{&&\hat{y}\ar[ld]_{\epsilon}
\ar[rd]^{\xi}&&\\
&\hat{y}\ar[ld]_{\alpha}\ar[rrrd]^{\beta}&&\hat{y}\ar[llld]_{\op{id}_{\hat{y}}}\ar[rd]^{\op{id}_{\hat{y}}}&\\
\hat{y}&&&&\hat{y}\,.}\een 
It is possible to show that the inverse for $(\op{id}_{\hat{x}},\Sigma)$ does not depend of the morphism $\delta_{\hat{y}}$.
The universal properties of the category of fractions and the quotient category imply the proposition. 
\end{proof}

A useful result for the relations \eqref{identification} is the following.
\begin{prop}\label{inv22}
For $\Sigma$ a $G$-cobordism with outgoing boundary the empty $G$-bundle, the relation $\op{id}_{\hat{x}}\sqcup (\Sigma\circ\overline{\Sigma})\sim\overline{\Sigma}\circ\Sigma$ is implied from the relation $\op{id}_{1}\sqcup (D\circ\overline{D})\sim\overline{D}\circ D$ for $D$ the disc $D:1\rightarrow 0$.
\end{prop}
\begin{proof}
It enough to assume $\Sigma$ a $G$-cobordism with connected base space from $x\in G$ to the empty $G$-bundle.
There is a unique $G$-cobordism $\Sigma':x\rightarrow 1$ such that $\Sigma=D\circ \Sigma'$. Consequently, we obtain 
\begin{equation}
    \overline{\Sigma}\circ\Sigma  =\overline{\Sigma'}\circ\overline{D}\circ D\circ\Sigma' =\overline{\Sigma'}\circ\op{id}_1\circ\Sigma'\sqcup S
    =\overline{\Sigma}\sharp\Sigma\sqcup S
\end{equation}
where $\overline{\Sigma}\sharp\Sigma$ is the connected sum of $G$-cobordisms and $S=D\circ\overline{D}$ is the unique $G$-cobordism over the sphere.
The base space of $\overline{\Sigma}\sharp\Sigma$ is a connected surface with exactly two boundary circles. For this space, we can find a simple closed curve which separates the space in two parts where one part contains the two boundary circles and it is maximal in the sense that the induced $G$-bundle is trivial.
In fact, after we fill with discs the two boundary circles produced by the curve, we obtain two disconnected $G$-cobordisms where one is the identity $G$-cylinder $\op{id}_x$ and the other is the $G$-cobordism $\Sigma\circ\overline{\Sigma}$. Notice that we have used for a second time the relation $\op{id}_{1}\sqcup (D\circ\overline{D})\sim\overline{D}\circ D$. This shows the proposition.
\end{proof}

\begin{thm}\label{confun}$\pi_1(B\cobG)=\ZZ \oplus \Omega_2^{SO}(BG)$. 
\end{thm}
\begin{proof}
For $\Gru$ the groupoid of fractions of $\cobG$, the fundamental group $\pi_1(B\cobG)$ can be identified with the automorphism group $\op{Hom}_{\Gru}(0,0)$, where $0$ denotes the trivial $G$-bundle. 
This monoid is composed by roofs $(\Gamma,\Sigma)$ where $\Gamma$ and $\Sigma$ are $G$-cobordisms with closed base space, inside the quotient category $\cobG_\sim$. The composition of two roofs $(\Gamma_1,\Sigma_1)$ and $(\Gamma_2,\Sigma_2)$ is given by $(\Gamma_1\sqcup \Gamma_2,\Sigma_1\sqcup\Sigma_2)$ which notice is commutative. Now we show there is the following short exact sequence 
\ben\label{guero}
0\rightarrow \ZZ\rightarrow\pi_1(B\cobG)\rightarrow \Omega_2^{SO}(BG)\rightarrow 0\,.
\een
We define the homomorphism $\pi_1(B\cobG)\rightarrow \Omega_2^{SO}(BG)$ which for a roof $(\Gamma,\Sigma)$ we assign the equivariant bordism class 
$[\Sigma]-[\Gamma]$. In order to show this assignment is well defined we need that the relation \eqref{identification} for the quotient category $\cobG_\sim$, is invariant up to equivariant bordism, however, this is precisely the Proposition \ref{inv22}. In fact, the relation  
$\op{id}_{1}\sqcup (D\circ\overline{D})\sim\overline{D}\circ D$ turns into the property that the operations of connected sum and disjoint union are equivalent in equivariant bordism, only that there is a fiber produced by the $G$-cobordism over the sphere.
The monomorphism $\ZZ\rightarrow\pi_1(B\cobG)$ is induced by $1\Tot (\op{id}_0,S)$ where $S$ is the $G$-cobordism over the sphere which is trivial in $\Omega_2^{SO}(BG)$.

In order to show the splitting of \eqref{guero}, we consider the map $\Omega_2^{SO}(BG)\rightarrow \pi_1(B\cobG)$ which starts with a  $G$-cobordism with closed base space and we kill all the trivial monodromy. This means that for every curve in the base space with trivial monodromy, we separate along this curve by capping with two discs. This extends to the total space by a trivial bundle. In \cite[pag. 7]{Domi} the author together with Dom\'inguez, show that the group $\Omega_2^{SO}(BG)$ can be obtained by considering the $G$-cobordisms with closed base space up to kill trivial monodromy. In addition, this map keeps in mind the number of times we separate by a trivial monodromy, by adding the negative of the $G$-cobordism over the sphere. Notice this process of separating and adding the negative of the $G$-cobordism over the sphere, is the inverse of the connected sum identification from previous paragraph. Thus the map $\Omega_2^{SO}(BG)\rightarrow \pi_1(B\cobG)$ is an splitting of the exact sequence \eqref{guero} and the proposition follows.
\end{proof}

Finally, the $G$-equivariant version of the embedded cobordism category of \cite{a6}, denoted by $\mathcal{C}_2^G$, has connected components $\pi_0(B\mathcal{C}_2^G)\cong \Omega_1^{SO}(BG)=G/[G,G]$. B\"okstedt-Svane \cite{a666} show that 
the fundamental group $\pi_1(B\mathcal{C}_2^G)$ is obtained in terms of generators up to the ``chimera relations". 
We deduce in the next paragraph these relations using the identification \eqref{identification}. This means that the canonical functor $\mathcal{C}_2^G\longrightarrow\cobG$ induces an isomorphism between the fundamental groups.

Take $\Sigma_i$, $i\in\{1,2,3,4\}$, $G$-cobordisms from $X$ to the empty $G$-bundle. We obtain the following relations for roofs:
\begin{align*}
    (\op{id}_0,\Sigma_3\circ\overline{\Sigma_1}\sqcup\Sigma_4\circ\overline{\Sigma_2}) &  = \left(\Sigma_1\circ\overline{\Sigma_1}\sqcup\Sigma_2\circ\overline{\Sigma_2}, (\Sigma_3\circ \overline{\Sigma_2}\circ \Sigma_2\circ\overline{\Sigma_1})\sqcup(
    \Sigma_4\circ \overline{\Sigma_1}\circ \Sigma_1\circ\overline{\Sigma_2}
    )\right)  \\
    &
    = \left(\Sigma_1\circ\overline{\Sigma_1}\sqcup\Sigma_2\circ\overline{\Sigma_2}, (\Sigma_3\circ \overline{\Sigma_2}\sqcup \Sigma_2\circ\overline{\Sigma_1})\sqcup(
    \Sigma_4\circ \overline{\Sigma_1}\sqcup \Sigma_1\circ\overline{\Sigma_2}
    )\right)\\
    &=  \left(\Sigma_1\circ\overline{\Sigma_1}\sqcup\Sigma_2\circ\overline{\Sigma_2}, (\Sigma_3\circ \overline{\Sigma_2}\sqcup \Sigma_4\circ\overline{\Sigma_1})\sqcup(
    \Sigma_2\circ \overline{\Sigma_1}\circ \Sigma_1\circ\overline{\Sigma_2}
    )\right) \\&
    = \left(\op{id}_0, \Sigma_3\circ \overline{\Sigma_2}\sqcup \Sigma_4\circ\overline{\Sigma_1}
    \right)\,.
\end{align*}

\section{Analysis of subcategories}
\label{sec2}

In this section we study the classifying space of the following subcategories of $\cobG$:
\begin{enumerate}
\item $\cobG_0$ is the full subcategory of $\cobG$ with only one object given by the empty $G$-bundle.
\item $\cobG_{>0}$ is the subcategory of $\cobG$ with the same objects of $\cobG$ except for the empty $G$-bundle and where for the base space,
each connected component of every morphism has non empty incoming boundary and non empty outgoing boundary.
\item $\cobG_b$ is the subcategory of $\cobG$ with the same objects of $\cobG$ and where for the base space,
each connected component of every morphism has non empty outgoing boundary.
\item $\cobG_1$ is the full subcategory of $\cobG_{>0}$ where 
each object is a principal $G$-bundle with base space the circle.

 \end{enumerate}

We denote by $\cob$, $\cob_{>0}$, $\cob_b$ and $\cob_1$ the cobordism categories associated to the trivial group.
Tillmann \cite{a3} shows the classifying space of $\cob_0$ is the infinite torus $T^\infty$. The same is true for $\cobG_0$.

\begin{prop}\label{proposition1} The classifying space $B\cobG_0$ is homotopic to the infinite dimensional torus $T^\infty$.
\end{prop}

\begin{proof}
This is similar to the proof of Tillmann \cite{a3}. Denote by the sequence $s_0$, $s_1$, $s_2$, $\cdots$ the $G$-cobordism classes with base space the sphere, the tori, a closed surface of genus 2 and so on. Since the group $G$ is finite, this is a sequence of finite positive integer numbers. Thus $\cobG_0$ is isomorphic to the direct sum $\bigoplus_{i\geq 0} \NN^{s_i}$ 
and since the classifying space of $\NN$ is the circle, then the classifying space $B\cobG_0$ is the infinite torus.
\end{proof}

Tillmann \cite{a3} defines a functor $\Phi:\cob_{>0}\rightarrow\cob_1$, which is the constant map in objects and each morphism $\Sigma$ with $n$ incoming circles, $c$ connected components, genus $g$ and $m$ outgoing circles, is mapped by $\Phi$ to
\ben
\Phi(\Sigma):=\frac{1}{2}(m-n-\chi(\Sigma))=g+m-c\,.
\een
This represents the unique morphism in $\cob_1$ with genus $g+m-c$. 
Moreover, $\Phi$ is left adjoint to the inclusion $\cob_1\hookrightarrow\cob_{>0}$ with a natural transformation defined for the object $n$, by the morphism $P_n:n\rightarrow 1$ which is the pair of pants with $n$ legs. Thus we have the commutativity of
\ben\label{pasa}
\xymatrix{n\ar[r]^{P_n}\ar[d]_{\Sigma} & 1\ar[d]^{g+m-c}\\ m \ar[r]^{P_m} & 1}
\een
Furthermore, the functor $\Phi$ can be extended to $\cob_b$. Since the subcategory $\cob_1$ is isomorphic to the natural numbers $\NN$
(the category with a single object and hom-set the natural numbers)
, the classifying spaces of $\cob_{>0}$, $\cob_b$ and $\cob_1$ have the same homotopy type of the circle $S^1$.

\begin{rem}\label{remi}
In order to implement the arguments of Tillman \cite{a3} for a finite group $G$, we need the following assumption:
every $G$-cobordism with connected base space and non-empty incoming boundary $\hat{x}=(x_1,\cdots,x_n)$, with $x_i\neq 0$ for $1\leq i\leq n$, factorises through the precomposition of the $G$-cobordism $P_{\hat{x}}$ (the pair of pants with multiple legs).
It would help to see an example in Figure \ref{f1}. This assumption will be taken until the end of the paper.

\begin{figure}
    \centering

   \includegraphics{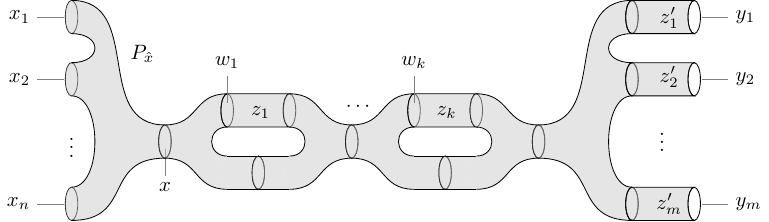}
    
    \caption{A morphism in $\cobG$.}
    \label{f1}
\end{figure}

\end{rem}



For a finite group $G$, we define the functor $\Phi^G:\cobG_{>0}\rightarrow\cobG_1$ in objects 
$(x_1,\cdots,x_n)$ by the product $x=\prod_{i=1}^n x_i$. For the morphisms of $\cobG_{>0}$, we consider first the case where we have a connected base space. Take $\Sigma$ a $G$-cobordism (with connected base space) which goes from $\hat{x}:=(x_1,\cdots,x_n)$ to $\hat{y}:=(y_1,\cdots,y_m)$. By hypothesis, see Definition \ref{cobordismoe}, the morphism $\Sigma$ factors by a precomposition of the $G$-cobordism over a pair of pants $P_{\hat{x}}$. 
Therefore, 
there is a unique $G$-cobordism $\Phi^G(\Sigma)$, from $x=\prod_{i=1}^nx_i$ to $y=\prod_{j=1}^my_j$, 
in the subcategory $\cobG_1$, satisfying the following commutative diagram: 
\ben \label{dfin}
\xymatrix{\hat{x} \ar[r]^{P_{\hat{x}}}\ar[d]_{\Sigma}&x\ar[d]^{\Phi^G(\Sigma)}\\ \hat{y}\ar[r]^{P_{\hat{y}}} & y}
\een
\begin{figure}
    \centering

    \includegraphics{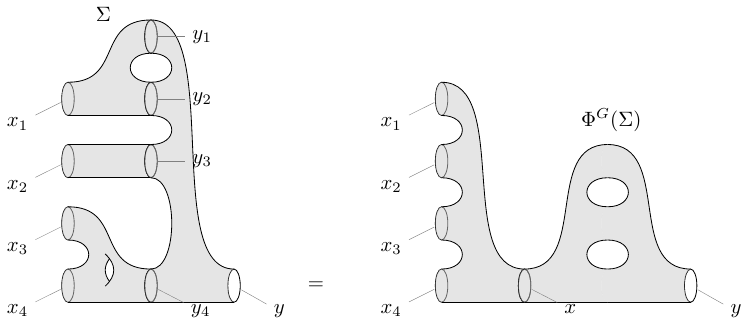}
  
    \caption{For $\Sigma:(x_1,x_2,x_3,x_4)\rightarrow(y_1,y_2,y_3,y_4)$ we represent the identity $P_{(y_1,y_2,y_3,y_4)}\circ \Sigma= \Phi^G(\Sigma)\circ P_{(x_1,x_2,x_3,x_4)}$.}
    \label{figdiaf}
\end{figure}
A monoidal property (without unit) is obtained for morphisms $\Sigma_1,\cdots,\Sigma_r$, where $\Sigma_i$, $1\leq i\leq r$, is a $G$-cobordism in $\cobG_{>0}$ with connected base space, which goes from $\hat{x}_i:=(x_{i1},\cdots,x_{in_i})$ to $\hat{y}_j:=(y_{j1},\cdots,y_{jm_j})$ ($x_i=\prod_{k=1}^{n_i}x_{ik}$ and $y_j=\prod_{k=1}^{m_j}y_{jk}$). Therefore, there is a unique $G$-cobordism $\Phi^G(\Sigma_1\sqcup\cdots\sqcup\Sigma_r)$, from $x_1\cdots x_r$ to $y_1\cdots y_r$, 
in the subcategory $\cobG_1$, satisfying the following commutative diagram: 
\ben \label{dfin2}
\xymatrix{\hat{x}_1\sqcup\cdots\sqcup\hat{x}_r \ar[rr]^{P_{\hat{x}_1\sqcup\cdots\sqcup\hat{x}_r}}
\ar[d]_{\Sigma_1\sqcup\cdots\sqcup\Sigma_r}  & & x_1\cdots x_r\ar[d]^{\Phi^G(\Sigma_1\sqcup\cdots\sqcup\Sigma_r)}\\ 
\hat{y}_1\sqcup\cdots\sqcup\hat{y}_r
\ar[rr]^{P_{\hat{y}_1\sqcup\cdots\sqcup\hat{y}_r}} &  & y_1\cdots y_r}
\een
It is straightforward to see that $\Phi^G$ is a functor because of the uniqueness of the factorisations, see Figure \ref{figdiaf} for an example.
 As a consequence, we obtain the following theorem:

\begin{thm}\label{thm123} The inclusion functor $\cobG_1\hookrightarrow\cobG_{>0}$ has a left adjoint $\Phi^G$.
\end{thm}

Likewise, we can extend the functor $\Phi^G$ to the category $\cobG_b$ and Theorem \ref{thm123} follows. Consequently, the categories $\cobG_{>0}$, $\cobG_b$ and $\cobG_1$ have homotopy equivalent  classifying spaces.

Finally, consider the $[G,G]$-graded monoid
$\ma{M}_G:=\{\ma{M}_x\}_{x\in [G,G]}$,
where $\ma{M}_x$ is the collection of $G$-cobordisms over handlebodies with monodromy $x\in G$ over the boundary circle. There are products $\ma{M}_{x_1}\times \ma{M}_{x_2}\longrightarrow\ma{M}_{x_1x_2}$ given by the composition with the $G$-cobordism over the pair of pants  $P_{(x_1,x_2)}$. 
Also, there is an action $\rho_z:\ma{M}_x\longrightarrow \ma{M}_{zxz^{-1}}$, $z\in G$, defined by the composition with a $G$-cobordism over the cylinder with entry $x$ and conjugation by $z$. The product is associative (with unit the disc) and twisted commutative and indeed, we mimic the axioms of a $G$-Frobenius algebra \cite{a10,GS}. 

The importance of the monoid $\ma{M}_G$ is that it works as a model for the composition of $\cobG_1$. This is because every morphism in $\cobG_1$ can be written as the gluing of two $G$-cobordisms, one over a handlebody and the other over a pair of pants. Thus the composition of two elements in $\cobG_1$ is given by the product of their associated $G$-cobordisms over the handlebodies.

An important part is the submonoid associated to the neutral element, given by $\ma{M}_1$. This consists of $G$-cobordisms over handlebodies with trivial monodromy over the boundary circle. 
This monoid is abelian, finitely generated and without torsion. 
Consequently, $\ma{M}_1$ is a finite direct sum $\NN^{r(G)}$ for some positive integer $r(G)$.
This number $r(G)$ writes as a finite sum $r_1(G)+r_2(G)+\cdots$, where $r_i(G)$ is the number of generators with base space of genus $i$. For instance, $r_1(G)$ coincides with the cardinality of the quotient of the commuting pairs by an action of the special linear group. When $G$ is an abelian finite group we obtain $r(G)=r_1(G)$, this is a multiplicative function and to find the exact value is feasible. We conclude the section with the following result.

Denote by $\NN^{r(G)}$ the category with a single object and hom-set the disjoint union of $r(G)$-copies of the natural numbers.

\begin{prop}\label{popis}
For $G$ a finite abelian group, the subcategory $\cobG_1$ is isomorphic with a finite direct sum $\NN^{r(G)}$.
\end{prop}

\section{Multiplicative structure}
\label{sec4}

In this section we consider the category $\cobG$ and the subcategory $\cobG_{>0}$ from Section \ref{sec2} with the assumption given in the Remark \ref{remi}. In addition,
these categories are restricted to the full subcategories where the objects are 
in the connected component of the empty $G$-bundle (or the trivial $G$-bundle over the circle). We denote these subcategories by the same symbols $\cobG$ and $\cobG_{>0}$.

\begin{figure}
    \centering
\includegraphics{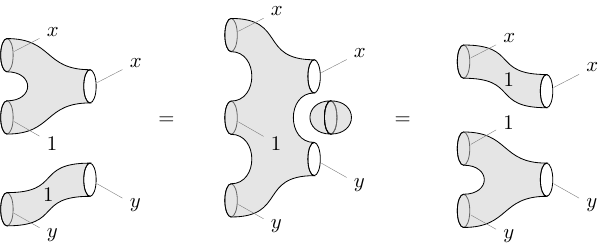}
    \caption{Unit axiom for the symmetric monoidal category $\tilde{\cobG_{>0}}$.}
    \label{unit}
\end{figure}

The category $\cobG_{>0}$ does not have a unit for the monoidal structure.  
However, there is a way to get a unit using the trivial $G$-bundle over the sphere and constructing a quotient category $\tilde{\cobG_{>0}}$ as follows: 
first of all,  principal $G$-bundles with connected base space admit a connected sum which is a $G$-equivariant extension of the usual connected sum of surfaces. 
For the product $\cobG_{>0}\times \ZZ$,
take morphisms $\gamma_i$, $i\in \{1,2\}$, in $\cobG_{>0}$, each one with connected base space, with the relation $(\gamma_1\sqcup\gamma_2,n)\sim (\gamma_1\sharp\gamma_2,n+1)$, where $\gamma_1\sharp\gamma_2$ denotes the connected sum of the morphisms $\gamma_1$ and $\gamma_2$.
Consider the equivalence relation in $\cobG_{>0}\times \ZZ$ generated by these identifications using 
the composition and the monoidal structure and denote the quotient category by 
$\tilde{\cobG_{>0}}$. In Figure \ref{unit}, we illustrate the triangle diagram for the unit axiom. 
Thus $\tilde{\cobG_{>0}}$ is a symmetric monoidal category.

We have a symmetric monoidal functor $\Psi:\cobG \longrightarrow\tilde{\cobG_{>0}}$, which for objects $(x_1,\cdots, x_n)$ with $x_i\neq 0$, $1\leq i\leq n$, the functor is the identity. As one might expect the empty $G$-bundle is sent to the trivial $G$-bundle over the circle, i.e., $0\stackrel{\Psi}{\longrightarrow}1$. 
Denote by $S$ the unique $G$-cobordism over the sphere.
 For the morphisms, 
 $\Psi$ is defined first for the disc $D:1\rightarrow 0$ and the reverse disc $\overline{D}$, by $\Psi(D)=(\op{id}_1,1)$ and $\Psi(\overline{D})=(\op{id}_1,0)$, where $\op{id}_1$ is the identity over the trivial $G$-bundle over the circle. 
  Now, for $\Sigma$ a morphism in $\cobG$ from $(x_1,\cdots,x_n)$ to $(y_1,\cdots,y_m)$, with $x_i,y_j\neq 0$, $i\in\{1,\cdots,n\}$ and $j\in\{1,\cdots,m\}$, we consider that the base space of $\Sigma$ has $k$ connected components, so we can write $\Sigma=\Sigma_1\sqcup \cdots\sqcup \Sigma_k$, where $\Sigma_l$, $l\in\{1,\cdots,k\}$, is the $G$-cobordism over the $l$ connected component. Thus we take $\Psi(\Sigma)=(\Sigma_1\sharp\cdots\sharp\Sigma_k,k-1)$. Using the composition of $\cobG$ we have a well defined functor $\Psi$ and it is not so difficult to see that the functor is symmetric and monoidal.

\begin{prop}
There is a symmetric monoidal functor $\Psi:\cobG \longrightarrow\tilde{\cobG_{>0}}$.
\end{prop}

Similarly as in Section \ref{sec2}, the subcategory $\cobG_1\subset \cobG_{>0}$ defines the subcategory $\tilde{\cobG_{1}}$. Again, the inclusion 
$\tilde{\cobG_{1}}\rightarrow\tilde{\cobG_{>0}}$ has a left adjoint denoted by 
\ben
\tilde{\Phi}^G:\tilde{\cobG_{>0}}\longrightarrow\tilde{\cobG_{1}}\,.
\een
The natural transformation is provided by the $G$-cobordism over the pair of pants with multiple legs. Thus for $\Sigma:\hat{x}\rightarrow\hat{y}$ a morphism in 
$\tilde{\cobG_{>0}}$, we obtain the following commutative square
\ben \label{dfina}
\xymatrix{\hat{x} \ar[r]^{P_{\hat{x}}}\ar[d]_{\Sigma}&x\ar[d]^{\tilde{\Phi}^G(\Sigma)}\\ \hat{y}\ar[r]^{P_{\hat{y}}} & y}
\een

The extension of the functor $\Phi^G$ to the whole category $\cobG$ is defined by the following composition:
\ben\label{extension}
\xymatrix{\cobG \ar[r]^{\Psi}\ar[rd]_{{{\Phi}'^G}}& \tilde{\cobG_{>0}}\ar[d]^{\tilde{\Phi}^G}\\ & \tilde{\cobG_{1}}}
\een
Now, we know $\Psi$ and $\tilde{\Phi}^G$ are symmetric monoidal functors and hence the induced maps on classifying spaces are maps of infinite loop spaces. However the composition 
\ben
\cobG_1\hookrightarrow \cobG\stackrel{{\Phi'}^G}{\longrightarrow}\tilde{\cobG_1}
\een
is almost never a homotopy equivalence in classifying spaces. For example, in the case $G$ is an abelian group we can show that the group completion of $\tilde{\cobG_1}$ is $\ZZ\oplus \Omega_2^{SO}(BG)$ while $\cobG_1$ is the finite sum $\NN^{r(G)}$. Therefore, in the equivariant case the argument of Tillmann \cite{a3} can not be entirely 
applied to produce an splitting.

This is what we can say for the classifying space of whole category $\cobG$.

\begin{thm}Every connected component of the $G$-cobordism category has the homotopy type of the product $S^1\times X^G$ where $\pi_1(X^G)=\Omega_2^{SO}(BG)$.
\end{thm}
\begin{proof}
Consider the subcategory of $\cobG_1$ where every element is represented by a trivial bundle. Notice this subcategory is isomorphic with the natural numbers $\NN$. Because of the definition of $\tilde{\cobG_1}$ we have a projection to $\ZZ$ and the composition \ben\NN\hookrightarrow\cobG\stackrel{{\Phi'}^G}{\longrightarrow}\tilde{\cobG_1}\rightarrow\ZZ\een is just the inclusion $\NN\hookrightarrow\ZZ$ which is a homotopy equivalence in classifying spaces. 
\end{proof}


\section{Applications to $G$-TQFT's}
\label{sec5}

First, consider functors $F:\cobG\longrightarrow\ma{D}$ with $\ma{D}$ a groupoid. The universal property of the groupoid of fractions $\cobG[{\cobG}^{-1}]$ 
associates a unique functor 
$\overline{F}:\cobG[{\cobG}^{-1}]\longrightarrow \ma{D}$, satisfying the commutative diagram 
\ben
\xymatrix{
\cobG \ar[rrd]_F\ar[rr]^{P_{\cobG}}& & \cobG[{\cobG}^{-1}]\ar[d]^{\overline{F}}\\ && \mathscr{D}}
\een
Similarly, the groupoid $\cobG[{\cobG}^{-1}]$ can be considered as the union of its connected components given by $G/[G,G]$. Any two connected components are related by the $G$-cobordisms over the pair of pants. Therefore, $F$ is determined by the restriction to some connected component plus the image of the $G$-cobordisms over the pair of pants.
Furthermore, there is a factorisation of the localization functor $P_{\cobG}$ through the category $\tilde{\cobG_{>0}}$ from Section \ref{sec4}. Thus we have the following diagram 
\ben
\xymatrix{\cobG\ar[rr]^{P_{\cobG}}\ar[rd]_\Psi & & \cobG[{\cobG}^{-1}]\\ &\tilde{\cobG_{>0}}\ar[ru] &
}
\een
Denote by $T$ the functor $\tilde{\cobG_{>0}}\longrightarrow \cobG[{\cobG}^{-1}]$ and set $F'=\overline{F}\circ T$.
In Section \ref{sec4} we have an adjunction between the categories $\tilde{\cobG_{>0}}$ and $\tilde{\cobG_1}$, which was used to define the functor ${\Phi'}^G:\cobG\longrightarrow\tilde{\cobG_1}$. 
For $\Sigma$ a morphism in $\cobG$, we consider its image $\Psi(\Sigma):\hat{x}\longrightarrow\hat{y}$, with $x=\prod x_i$ and $y=\prod y_j$. This adjunction has the following commutative diagram 
\ben
\xymatrix{\hat{x} \ar[d]_{\Psi(\Sigma)}\ar[r]^{P_{\hat{x}}}& x\ar[d]^{{\Phi'}^G(\Sigma)}\\ \hat{y}\ar[r]^{P_{\hat{y}}} & y}
\een
Therefore, the functor $F:\cobG\longrightarrow\mathscr{D}$ is completely determined as follows:
\ben
F(\Sigma)=F'(P_{\hat{y}})^{-1}\circ F'({\Phi'}^G(\Sigma))\circ F'(P_{\hat{x}})\,.
\een

As an application we consider an abelian group $\ma{D}$. We recover the case of Tillmann \cite{a3} by setting $a_0=F(S^2)$ and $c_n=F(P_n)$. Thus for $\Sigma:n\rightarrow m$ the functor $F$ is given by 
\ben
F(\Sigma)=c_n-c_m+a_0(\Phi(\Sigma))\,,
\een
where $b_n+c_n=b_1$ with $b_n=F(D_n)$ for $D_n:0\rightarrow n$ the disjoint union of $n$ disks.

Second, we consider symmetric monoidal functors $F:\cobG\rightarrow\op{Vect}_\CC$, where $\op{Vect}_\CC$ is the category of finite dimensional complex vector spaces. 
They are currently called $G$-topological quantum field theories \cite{a10}. 
A folk result states a correspondence between every symmetric monoidal functor $F:\cobG\rightarrow\op{Vect}_\CC$ and a $G$-Frobenius algebra \cite{a10,GS}. 
Briefly, a $G$-Frobenius algebra is a $G$-graded vector space $A=\bigoplus_{x\in G} A_x$, where $A_x$ is a finite dimensional $\CC$-vector space, with $G$-graded associative products $A_x\otimes A_y\rightarrow A_{xy}$ which are twisted commutative, also there is a trace $\theta:A_1\rightarrow\CC$ which is non degenerate and we have an action $\rho:G\rightarrow\op{Aut}(A)$ of algebra automorphisms and some other axioms.  
We are interested when $F$ is morphism inverting, hence  $A_x\cong \CC$ for any $x\in G$. The algebraic structure is called discrete torsion determined by Turaev \cite{turaev} as follows:
\ben
\CC_b(G):=\bigoplus_{x\in G}\CC\times\{x\}\,,
\een
where for a basis $e_x$ in each $A_x$ and $b\in H^2(G,\CC^*)$ a 2-cocycle, the product is $e_x\cdot e_y=b(x,y)\cdot e_{xy}$.
There is a commutative diagram 
\ben
\xymatrix{
\cobG \ar[rrd]_F\ar[rr]^{P_{\cobG}}& & \cobG[{\cobG}^{-1}]\ar[d]^{\overline{F}}\\ && \op{Vect}_\CC}
\een
Thus $F$ is completely determined by the restriction to the automorphism group of
the empty $G$-bundle in 
$\cobG[{\cobG}^{-1}]$, plus the images of a generating set of the monoid $\ma{M}_G$ from Section \ref{sec2}. This induces a representation of the fundamental group $\ZZ\oplus \Omega_2^{SO}(BG)\rightarrow \CC^*$ from Section \ref{sec3}.  

\begin{figure}
    \centering

    \includegraphics{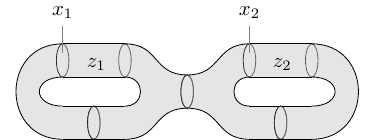}
    
        \caption{Example of a $G$-cobordism over a closed surface.}
    \label{corte}
\end{figure}

Finally, for any $G$-topological quantum field theory $F:\cobG\rightarrow\op{Vect}_\CC$, we consider the subset of morphisms $\tilde{\cobG_0}$ used in Section \ref{sec3}. This consists of those endomorphisms in $\cobG$ which are the disjoint union of identity $G$-cobordisms over cylinders and $G$-cobordisms over closed surfaces. 
Similarly, as in the case of Tillman \cite{a3}, any $G$-cobordism over a closed surface can be written as a double construction given by a $G$-cobordism over the gluing $M\sqcup_X \overline{M}$, where 
$M$ is a surface with boundary $X$ and $\overline{M}$ is $M$ with opposite orientation. The manifold $M$ is found observing that any $G$-cobordism over a closed surface has a decomposition in terms of $G$-cobordisms over pairs of pants and cylinders, see Figure \ref{corte} for an example.
Then we cut the surface in two halves given by the front and the back side observing they have the same $G$-cobordisms but with different orientation.
This implies that any $G$-cobordism over a closed surface is assigned to a non-zero complex number by the non-degeneracy of the inner product. 
Consequently, we obtain that every $G$-topological quantum field theory $F:\cobG\rightarrow\op{Vect}_\CC$ factors through $\cobG[{\tilde{\cobG_0}}^{-1}]$.


\newpage
{\bf Acknowledgement}

The author is supported by c\'atedras CONACYT and Proyecto CONACYT ciencias b\'asicas 2016, No. 284621.

\bibliographystyle{amsalpha}
\bibliography{biblio}

\end{document}